\documentclass[bezier,11pt]{article}
\usepackage{amsfonts}
\usepackage{amsmath, amssymb}
\usepackage{mathrsfs}
\usepackage{rotating}
\usepackage{amsthm} 
\usepackage{bbm}

\newtheorem{theorem}{Theorem}[section]
\newtheorem{corollary}[theorem]{Corollary}
\newtheorem{assumption}[theorem]{Assumption}

\newtheorem{lemma}[theorem]{Lemma}

\newtheorem{example}[theorem]{\bf Example}

\numberwithin{equation}{section}

\newcommand{\lk}{ [\kern-0.1em[  }
\newcommand{\rk}{ ]\kern-0.1em] }
\newcommand{\0}{\cdot }
\newcommand{\lag}{\langle }
\newcommand{\rag}{\rangle }

\newcommand{\dint}{\displaystyle\int}

\newcommand{\E}{{\cal E}}

\newcommand{\ttau}{{\widetilde\tau}}

\newcommand{\wt}{\widetilde}
\newcommand{\wh}{\widehat}
\newcommand{\C}{{\widetilde C}}

\newcommand{\Q}{{\!{^Q}}}

\newcommand{\p}{{^{(p)}}}

\newcommand{\scr}{\mathscr }
\newcommand{\bb}{\mathbb }

\begin{document}

\def \R {{\bf R}}
\def \C {{\bf I}}
\def \E {{\bf E}}
\def \N {{\bf N}}
\def \P {{\bf P}}
\def \Q {{\bf Q}}
\def \bb {\mathbb}
\def \l{\langle}
\def \r{\rangle}
\def \p{\partial}
\def \ba{\begin{eqnarray*}}
\def \ea{\end{eqnarray*}}
\def \bc{\begin{center}}
\def \ec{\end{center}}
\def \bl{\begin{flushleft}}
\def \el{\end{flushleft}}
\def \br{\begin{flushright}}
\def \er{\end{flushright}}
\def \bd{\begin{document}}
\def \ed{\end{document}}

\title{ \LARGE Parametrization in the progressively enlarged filtration
\thanks{\noindent Supported by National Natural Science Foundation of China(no.11171215), National Nature Science
of Shanghai(no.13ZR1422000)
and Shanghai 085 Project. }}
\author{ Kun Tian$^{1}$ \thanks{\noindent The
corresponding author, his Email is tkrp1980@sjtu.edu.cn.},\quad Dewen Xiong $^{1}$
\quad and\quad Zhongxing Ye $^{1,2}$ \\
\noindent {\small 1. Department of Mathematics, Shanghai Jiao Tong
University, Shanghai 200240, }
\\
\noindent {\small    China}
\\{\small 2. School of Business Information Management, Shanghai University of International}
\\{\small Business and Economics, Shanghai
201620, China}
}
\date{\,}
\maketitle \thispagestyle{empty}

\begin{abstract}
In this paper,
we assume that the filtration $\bb F$ is
generated by a $d$-dimensional Brownian motion
$W=(W_1,\cdots,W_d)'$ as well as an integer-valued random measure
$\mu(du,dy)$. The random variable $\ttau$ is the default time and $L$ is the default loss.
Let $\mathbb G=\{\scr
G_t;t\geq 0\}$ be the progressive enlargement of $\bb F$ by $(\ttau,L)$, i.e,
$\bb G$ is
the smallest filtration including $\bb F$ such that $\ttau$ is a $\bb G$-stopping
time and $L$ is $\scr G_\ttau$-measurable. We parameterize the conditional density process, which allows us  to describe the survival process $G$ explicitly.  We also obtain the explicit $\bb
G$-decomposition of a $\bb F$ martingale and  the predictable representation theorem  for a $(P,\bb
G)$-martingale by all known parameters. Formula parametrization in the enlarged filtration is a useful quality in financial modeling.

\end{abstract}

\section{Introduction}
In this paper, we assume that the filtration $\bb F$ is generated by
a $d$-dimensional Brownian motion $W=(W_1,\cdots,W_d)'$ as well as
an integer-valued random measure $\mu(du,dy)$. It is the fact that once the default happens,
the default loss  immediately generated, we
let $\ttau\geq 0$ be
the default time and $L$ be the default loss, which are random
variables. We are interest in the enlarged filtration $\mathbb G=\{\scr G_t;t\geq 0\}$ which is the smallest
filtration including $\bb F$ such that $\ttau$ is a $\bb G$-stopping
time and $L$ is $\scr G_\ttau$-measurable, $\bb G$ is called the progressive
enlargement of $\bb F$ by $(\ttau,L)$. We can easily see that
$$\scr G_t=\bigcap_{u>t}\scr G_u^0,~~\text{where}~~\scr G_u^0=\scr F_u\vee\sigma(\ttau\wedge u)\vee\sigma(L\bb I_{\ttau\le
 u}).$$
The similar progressive enlargement of filtration can also be found in Dellacherie and Meyer(1978)$^{\text{\cite{Dellacherie-Meyer1978}}}$, Pham(2010)$^{\text{\cite{Pham-2010}}}$,
Kchia, Larsson and Protter(2011) $^{\text{\cite{Kchia-Larsson-Protter-01}}}$ $^{\text{\cite{Kchia-Larsson-Protter-02}}}$.
It is notable that the
progressive enlargement filtration $\bb G$ in this paper is slightly different from the
traditional progressive enlargement of filtration in the
literature, see e.g. in Jeulin(1980)$^{\text{\cite{Jeulin1980}}}$,
Jacod(1987)$^{\text{\cite{Jacod1987}}}$, Jeanblanc, Yor and
Chesney(2009) $^{\text{\cite{Jeanblanc-Yor-Chesney2009}}}$,
El Karoui, Jeanblanc and Jiao(2010)$^{\text{\cite{El
Karoui-Jeanblanc-Jiao2010}}}$,
 Jeanblanc and Le Cam(2009)$^{\text{\cite{Jeanblanc-Cam2009}}}$,
Jeanblanc and Song(2011)$^{\text{\cite{Jeanblanc-Song2011}}}$,
Callegaro, Jeanblanc and
Zargari(2010)$^{\text{\cite{Callegaro-Jeanblanc-Zargari(2010)}}}$
and Jeanblanc and Song(2012)$^{\text{\cite{Jeanblanc-Song2012}}}$,
among many others.
The $\bb G$-decomposition of a
$(P,\bb F)$ martingale and the representation of a $\bb G$-martingale  are the fundamental problems in the enlargement of filtration, which have been widely studied, see in \cite{Callegaro-Jeanblanc-Zargari(2010)}, \cite{Dellacherie-Meyer1978}, \cite{El Karoui-Jeanblanc-Jiao2010}, \cite{Jeanblanc-Cam2009}, \cite{Jeanblanc-Song2011},
\cite{Jeanblanc-Song2012}, and which
have many applications in mathematical finance, see Duffie and
Huang(1986)$^{\text{\cite{Duffie-Huang1986}}}$, Karatzas and
Pikovsky(1996)$^{\text{\cite{Pikovsky-Karatzas1996}}}$, Amendinger,
Becherer and
Schweizer(2003)$^{\text{\cite{Amendinger-Becherer-Schweizer2003}}}$,
 Ankirchner, Dereichner and
Imkeller(2005)$^{\text{\cite{Ankirchner-Dereichner-Imkeller2005}}}$,
Jiao and Pham(2009)$^{\text{\cite{Jiao-Pham2009}}}$ and Eyraud-Loisel(2010)$^{\text{\cite{Eyraud-Loisel2010}}}$, etc.

It is well known that for a general enlargement of filtration,
a $(P, \bb F)$-martingale might not be a $(P, \bb G)$-semimartingale.
here we adapt the stronger Jacod's hypothesis (see in Jacod(1987)$^{\text{\cite{Jacod1987}}}$, Amendinger(1999)$^{\text{\cite{Amendinger(1999)}}}$, Callegaro, Jeanblanc and
Zargari(2010)$^{\text{\cite{Callegaro-Jeanblanc-Zargari(2010)}}}$)
as the following:
\begin{assumption}\label{assumption-01}
i)  Let $\eta(ds,dl)$ be the law of $(\ttau,L)$ and the $\bb F$-regular conditional law of $(\ttau, L)$ is equivalent
to the law of $(\ttau, L)$, i.e
$$P(\ttau\in ds,L\in dl|\scr F_t)\sim \eta(ds,dl),\quad\text{for every }t\ge 0;$$
ii) $\eta(ds, dl)$ has no atoms.
\end{assumption}

From Assumption \ref{assumption-01} and similarly to  Jacod(1987)$^{\text{\cite{Jacod1987}}}$ or
Amendinger(1999)$^{\text{\cite{Amendinger(1999)}}}$, there also exists a
so-called
\textbf{conditional density}
$p_t(s, l)$ which is a $(P, \bb
F)$-martingale to describe $\bb F$-regular conditional law of $(\ttau,
L)$,
such that for every
$(s,l)\in\bb R^+\times\bb R$, $p(s,l)$ is a c\`{a}dl\`{a}g $(P,\bb
F)$-martingale and for any $B\in\scr B(\bb R^+\times\bb R)$,
$$P((\ttau,L)\in B|\scr F_t)=\int_{B}p_t(s,l)\eta(ds,dl),\quad \text{for every }t\ge 0, P\text{-a.s.}$$
then $p_0(s, l)=1$, for any $t\geq 0, l\in \bb R$.

Assuming 1.1 is a foundmental assumption, see
\cite{Callegaro-Jeanblanc-Zargari(2010)},
\cite{El
Karoui-Jeanblanc-Jiao2010} and more recently $\cite{Kchia-Larsson-Protter-02}$.
If $m$ is a $(P, \bb F)$ martingale, we know from \cite{Kchia-Larsson-Protter-02} that $m^\bb G_t:=m_t-\dint_0^{t\wedge\ttau}\frac{d\lag m, \mu\rag_s+dJ_s}{G_{s-}}\bb I_{\ttau>t}-\dint_{t\wedge\ttau}^tdA_s$
is a $(P, \bb G)$-martingale (see Theorem 3 in
 \cite{Kchia-Larsson-Protter-02}), where $\mu$ is the martingale part of the conditional survive process $G_t=P(\ttau>t|\scr F_t)$ and $J$ is the dual predictable projection of $\Delta M_\ttau \bb I_{\lk\ttau,\infty\lk}$ onto $\bb G$,  from which one can obtain the decomposition of
$m$ in $\bb G$. However in practice, seeking $\mu$, $\lag m, \mu\rag$ and $A$ are very troublesome and it is hard to get this decomposition. Although this formula is  good, it is very inconvenient to calculate, which limits its application.

In this paper we will
investigate $\{p_t(s, l); t\geq
0\}$ more deeply, and rewrite the conditional density $p_t(s,l)$ into the following form
\begin{displaymath}
   \begin{aligned}
p_t(s, l)&=E[p_s(s,l)\big|\scr F_t]\bb I_{t<s}+p_s(s, l)\bb I_{t\geq s}
+\dint_0^tp_{u-}(s, l)\theta_1(u; s, l)'\bb I_{u>s}dW(u)\\
&+\dint_0^t\dint_Ep_{u-}(s, l)\theta_2(u, y; s, l) \bb
I_{u>s}\{\mu(du,dy)-\nu(du,dy)\},
 \end{aligned}
\end{displaymath}
where $\theta_{1}(\cdot;s,l)$ is an $\bb R^d$-valued $\bb F$-predictable process with $\theta_{1}(\cdot;s,l)\in L^2_{loc}(W)$ and $\theta_{2}(\cdot,\cdot;s, l)$
 is a $\widetilde{\scr P}(\bb F)$-measurable function with $\theta_{2}(\cdot,\cdot;s, l)\in G_{loc}(\mu) $.
We
show in Theorem \ref{theorem-p(s;s,l)} that
 $p_s(s, l)$, $\theta_1(u; s, l)$ and $\theta_2(u, y; s, l)$ must satisfy the following condition
\begin{equation}\label{condition-p(s;s,l)}
\begin{array}{ll}\dint_0^\infty\dint_\bb Rp_s(s,l)\bb
\eta(ds,dl)=1-\dint_0^\infty\Big\{\dint_0^u\dint_\bb
Rp_s(s,l)Z^{\theta_1,\theta_2}_{s,l}(u-)\theta_1(u;s,l)\eta(ds,dl)\Big\}'d
W(u)\\\qquad-\dint_0^\infty\dint_E\Big\{\dint_0^u\dint_\bb
Rp_s(s,l)Z^{\theta_1,\theta_2}_{s,l}(u-) \theta_2(u,y;s,l))
\eta(ds,dl)\Big\}\{\mu(du,dy)-\nu(du,dy)\}.\end{array}
\end{equation}
where
$$\begin{array}{ll}
Z^{\theta_1,\theta_2}_{s,l}(t):=\exp\bigg\{\dint_s^t\theta_1(u;s,l)'d
W(u)-\dfrac{1}{2}\dint_s^t\|\theta_1(u;s,l)\|^2du\bigg\}\\
\qquad\times\exp\bigg\{\dint_s^t\dint_E\big\{\ln(1+\theta_2(u,y;s,l))\big\}\mu(du,dy)
-\dint_s^t\dint_E\theta_2(u,y;s,l)\nu(du,dy)\bigg\}\end{array}$$
is a $(P, \bb F)$ martingale determined by $\theta_1$ and $\theta_2$. We view $p_s(s,l)$, $\theta_1$ and
$\theta_2$ as the parameters of the conditional density.
Inversely, given any parameters $\{p_s(s, l), \theta_1, \theta_2\}$ satisfying (\ref{condition-p(s;s,l)}), we can construct a new``conditional density", which may have
 applications in finance.

Based on the parameters $\{p_s(s, l), \theta_1, \theta_2\}$, we can give the  Doob-Meyer's decomposition of the survival process $G_t$ explicitly as the following
\begin{equation}\label{eqn-G-Doob-Meyer}
\begin{array}{ll}G_t
&=1-\dint_0^t\int_\bb Rp_s(s,l)\eta(ds,dl)\\
&\qquad-\dint_0^t\Big\{\dint_0^u\dint_\bb
Rp_s(s,l)Z^{\theta_1,\theta_2}_{s,l}(u-)\theta_1(u;s,l)\eta(ds,dl)\Big\}'d
W(u)\\
&\qquad-\dint_0^t\dint_E\Big\{\dint_0^u\dint_\bb
Rp_s(s,l)Z^{\theta_1,\theta_2}_{s,l}(u-) \theta_2(u,y;s,l))
\eta(ds,dl)\Big\}\{\mu(du,dy)-\nu(du,dy)\}.
\end{array}
\end{equation}
Then we consider the $\bb G$-decomposition of a $(P,\bb F)$
local martingale and give a
 $\bb G$-decomposition explicitly
depending on $p_s(s,l)$, $\theta_1(u;s,l)\bb I_{u>s}$ and
$\theta_2(u,y;s,l)\bb I_{u>s}$.  We also obtain the
predictable representation theorems both for a $(P,\bb
G)$-martingale and a $(P^*,\bb G)$-martingale, which reveals the relationship between the structure of $(P,\bb G)$-martingale and these parameters.

In short, the main results are:

1). An analysis of the conditional density process $p_t(s,l)$,
it is underlined that $p_t(s,l)$ is
completely determined by three parameters $p_t(s,l)$, $\theta_1$ and
$\theta_2$.

2). An explicit expression of the survival process $G_t$ in terms of the three parameters $p_t(s,l)$, $\theta_1$ and
$\theta_2$.

3). Theorem 3.4: The so-called enlargement filtration formula is proved. This is the formula
which gives the $\bb G$ semimartingale decompositions for the $\bb F$ local martingales.

4). Theorem 4.3, where the (weak) martingale representation formula in $\bb G$ is proved and formula parametrization is obtained.

\par The paper is organized as follows:
In Section 2, we characterize the conditional density process and
obtain a more explicit form of the Doob-Meyer's decomposition of the
survival process.
 In Section 3, we will first prove Lemma \ref{lemma-Af} and then
explicitly describe the $\bb G$-decomposition of a $(P, \bb
F)$-martingale. In  section 4, we obtain the martingale
representation theorem for  a $(P,\bb G)$-martingale. Conclusions are given in the end.

\section{The setup and notations}

We assume that $\bb F=\{\scr F_t; t\ge0\}$ is a filtration on $(\Omega,
\scr F, P)$ carrying an $d$-dimensional Brownian motion
$W=(W_1,\cdots,W_d)'$ as well as an integer-valued random measure
$\mu(du,dy)$ on $\bb R^+\times E$, where $(E,\scr E)$ is some auxiliary Blackwell space.

\begin{assumption}\label{assum-filtration-F}
The filtration $\bb F=\{\scr F_t; t\ge 0\}$ is the natural
filtration generated by $W$ and $\mu$, i.e.,
$$\scr F_t=\sigma\big\{W(s),\mu([0,s]\times A),B;
\quad 0\le s\le t,\; A\in\scr B,\; B\in\scr N\big\}$$ where $\scr N$
is the collection of $P$-null sets from $\scr  F$.
\end{assumption}

In the following, let $\scr P(\bb F)$ be the family of all
$\bb F$-predictable processes and $\wt{\scr P}(\bb F)=\scr P(\bb
F)\otimes\scr E$. We assume that the compensator of $\mu(du,dy)$ is
given by $\nu(du,dy)=F_u(dy)du$, where  $F_u(dy)$ is a transition
kernel from $(\Omega\times\bb R^+,\scr P(\bb F))$ to into $(E,\scr
E)$. It is known that $(y^2\wedge 1)*\nu\in \scr A^+_{loc} $, i.e, there exists a
sequence of increasing
stopping time $T_n$ with $\lim_{n\rightarrow \infty}T_n=\infty \ a.s.$, such that $E(\dint_0^{T_n}\dint_E(y^2\wedge 1)F_u(dy)du)<\infty$. Furthermore,  for any
local $(P,\bb F)$-martingale $M$ has the unique representation property
$$M_t=M_0+\dint_0^tf_1(u)'dW(u)+\dint_0^t\dint_Ef_2(u,y)(\mu(du,dy)-\nu(du,dy)),$$
where $f_{1}(u)$ is an $\bb R^d$-valued $\bb F$-predictable process with $f_1\in L^2_{loc}(W)$ and $f_{2}(u, y)$
 is a $\widetilde{\scr P}(\bb F)$-measurable function with $f_{2}\in G_{loc}(\mu) $ (see Lemma
 4.24 in page 185 of
Jacod and Shiryaev(1987) $^{\text{\cite{Jacod-Shiryaev1987}}}$ for more details).\\

\begin{remark}
For simply, here we only assume that $W$ is  an $d$-dimensional Brownian motion
In fact,  $W$ can be relaxed to any continuous local martingale.  As for integer-valued random measure
$\mu(du,dy)$,
the representation property of any local $(P,\bb F)$-martingale $M$ has
the more concise form without loss of generality.
\end{remark}

\begin{remark}
From  Proposition 2.1.14 and 2.1.15 in
 Jacod and Shiryaev(1987) $^{\text{\cite{Jacod-Shiryaev1987}}}$, one can see that there exists a $\bb F$-optional process
$\varphi=(\varphi_t)$ and a sequence of stopping times $(\hat\tau_k)$
such that for all positive $\wt{\scr P}(\bb F)$-measurable function
$W(\omega,t,y)$, $$\int_0^t\dint_E W(\omega,u,y)\mu(\omega; du, dy)
=\sum_{(k)} W(\hat \tau_k, \varphi_{\hat\tau_k} )\bb I_{\hat\tau_k\le
t}.$$ Furthermore, as the compensator of $\mu(du,dy)$ is
$\nu(du,dy)=F_u(dy)du$, so the filtration $\bb F$ is quasi-left
continuous.
\end{remark}

Let $\widetilde{\tau}$ be a non-negative random variable and $L$ be a
random variable on $(\Omega,\scr F)$, and let the $\mathbb G=\{\scr
G_t;t\geq 0\}$ be the smallest progressive enlargement filtration of $\bb F$ such that $\ttau$
is a $\bb G$-stopping time and $L$ is a $\scr G_\ttau$-measurable
random variable. Let $\eta(ds,dl)$ be the law of $(\ttau,L)$, i.e.,
$\eta(ds,dl)=P(\ttau\in ds,L \in dl)$, then $\dint_0^\infty\dint_\bb R\eta(ds,dl)=1$. We  adapt the stronger Jacod's hypothesis  as following
\\

\par\noindent \textbf{Assumption 1.1.}
i) the $\bb F$-regular conditional law of $(\ttau,L)$ is equivalent
to the law of $(\ttau,L)$, i.e.,
$$P(\ttau\in ds,L\in dl|\scr F_t)\sim \eta(ds,dl),\quad \text{for every }t\ge 0, P\text{-a.s.}$$
ii) $\eta(ds,dl)$ has no atoms.

Similar to  Proposition 1.3 in
Callegaro, Jeanblanc and Zargari(2010)$^{\text{\cite{Callegaro-Jeanblanc-Zargari(2010)}}}$, we have the following
lemma

\begin{lemma}
1) Any $\bb G$-predictable process $Y = (Y_t)_{t\ge0}$ is
represented as
$$Y_t=Y^0_t\bb I_{t\le\ttau}+Y^1_t(\ttau,L)\bb I_{t>\ttau},$$
where $Y^0$ is a $\bb F$-predictable process and where $Y^1_t(s,l)$
is a $\scr P(\bb F)\otimes\scr B(\bb R^+)\otimes\scr B(\bb R)$-measurable
function.\par
2) Any $\bb G$-optional process $Y = (Y_t)_{t\ge0}$ is
represented as
$$Y_t=Y^0_t\bb I_{t<\ttau}+Y^1_t(\ttau,L)\bb I_{t\ge\ttau},$$
where $Y^0$ is a $\bb F$-optional process and where $Y^1_t(s,l)$ is
a $\scr O(\bb F)\otimes\scr B(\bb R^+)\otimes\scr B(\bb R)$-measurable
function.
\end{lemma}

From Assumption \ref{assumption-01}, one can see from
Jacod(1987)$^{\text{\cite{Jacod1987}}}$ or
Amendinger(1999)$^{\text{\cite{Amendinger(1999)}}}$ that there
exists a strictly positive $\scr O(\bb F)\otimes\scr B(\bb R^+\times
\bb R)$-measurable function $(t;\omega; s,l)\longrightarrow
p_t(\omega; s,l)$, called the $(P,\bb F)$\textbf{-conditional
density} of $(\ttau,L)$ with respect to $\eta$, such that for every
$(s,l)\in\bb R^+\times\bb R$, $p(s,l)$ is a c\`{a}dl\`{a}g $(P,\bb
F)$-martingale and for any $B\in\scr B(\bb R^+\times\bb R)$,
$$P((\ttau,L)\in B|\scr F_t)=\int_{B}p_t(s,l)\eta(ds,dl),\quad \text{for every }t\ge 0, P\text{-a.s.}$$
then $p_0(s, l)=1$, for any $t\geq 0, l\in \bb R$.

By the
``change of probability measure" viewpoint of
 Song(1987)$^{\text{\cite{Song1987}}}$ and
similar to Callegaro, Jeanblanc and
Zargari(2010)$^{\text{\cite{Callegaro-Jeanblanc-Zargari(2010)}}}$,
we introduce the filtration $\bb G^{\ttau,L}=\{\scr G^{\ttau,L}_t;
t\ge 0\}$ with $\scr G^{\ttau,L}_t=\scr F_t\bigvee \sigma(\ttau,
L)$, one can see that $\bb G^{\ttau,L}$ is the initial enlargement
of the filtration $\bb F$ with $\ttau$ and $L$ and that $\bb
F\subset \bb G\subset \bb G^{\ttau,L}$. Let
$$Z_t=\dfrac{1}{p_t(\ttau,L)},$$
similar to
Grorud and Pontier(1998)$^{\text{\cite{Grorud-Pontier-1998}}}$  or
Amendinger(1999)$^{\text{\cite{Amendinger(1999)}}}$, one can see
that $Z$ is a strictly positive $(P,\bb G^{\ttau,L})$-martingale
with $E(Z_t)=1$, for every finite $t\ge 0$. Thus one can define a locally
equivalent probability measure $P^*$ by
$$\dfrac{dP^*}{dP}|_{\scr G^{\ttau,L}_t}=Z_t.$$

Similar to
Grorud-Pontier(1998)$^{\text{\cite{Grorud-Pontier-1998}}}$ or
Amendinger(1999)$^{\text{\cite{Amendinger(1999)}}}$, one can show
that
\begin{enumerate}
  \item under $P^*$, $(\ttau,L)$ is independent of $\scr F_t$ for
  every $t\ge 0$;
  \item $P^*|_{\scr F_t}=P|_{\scr F_t}$;
  \item $P^*|_{\sigma(\ttau,L)}=P|_{\sigma(\ttau,L)}$,
\end{enumerate}
which implies $P^*(\ttau\in ds,L\in dl|\scr F_t)=P(\ttau\in ds,L\in
dl)$. Similar to Lemma 1.4 in
Callegaro, Jeanblanc and Zargari(2010)$^{\text{\cite{Callegaro-Jeanblanc-Zargari(2010)}}}$,
we have the following lemma
\begin{lemma}\label{lemma-E(|Ft)}
1) Let $y_t(\ttau,L)$ be a $\scr G^{\ttau,L}_t$-measurable r.v.,
then for any $s\le t$,
$$E_{P^*}\big(y_t(\ttau,L)\big|\scr G^{\ttau,L}_s\big)=E_{P^*}\big(y_t(u,l)\big|\scr
F_s\big)\big|_{u=\ttau,\:l=L}\:;$$
2) if $y_t(\ttau,L)$ is
$P$-integrable, then
$$E\big(y_t(\ttau,L)\big|\scr G^{\ttau,L}_s\big)=\dfrac{1}{p_s(\ttau,L)}E\big(y_t(u,l)p_t(u,l)\big|\scr
F_s\big)\big|_{u=\ttau,\:l=L}\:.$$
\end{lemma}

We have the following Corollaries
\begin{corollary}
[Characterization of $(P, \bb G^{\ttau,L})$-martingales in terms of
$(P,\bb F)$-martingales] A process $y_t(\ttau,L)$ is a $(P,\bb
G^{\ttau,L})$-martingale if and only if $\{y_t(u, l)p_t(u, l); t\ge
0\}$ is a $(P,\bb F)$-martingale, for almost every $u\ge 0, l\in\bb
R$.
\end{corollary}

\begin{corollary}
Let $M=\{M_t;t\ge 0\}$ be a bounded $(P^*,\bb F)$-martingale, then
$M$ is a $(P^*,\bb G^{\ttau,L})$-martingale and hence a $(P^*,\bb
G)$-martingale.
\end{corollary}
\begin{proof}From part 1) of Lemma \ref{lemma-E(|Ft)}, one can see that for any $s\le t$
$$\begin{array}{ll}E_{P^*}\big(M_t\big|\scr G^{\ttau,L}_s\big)&=E_{P^*}\big(M_t\big|\scr
F_s\big)\big|_{u=\ttau,\:l=L}\\
&=M_s,\end{array}$$ thus $M$ is a $(P^*,\bb
G^{\ttau,L})$-martingale. Since $M_s$ is $\scr F_s$-measurable,
$M_s$ is $\scr G_s$-measurable, thus
$$\begin{array}{ll}E_{P^*}\big(M_t\big|\scr G_s\big)&=E_{P^*}\{E_{P^*}\big(M_t\big|\scr
G^{\ttau,L}_s\big)|\scr G_s\}\\
&=E_{P^*}\{M_s|\scr G_s\}\\
&=M_s,\end{array}$$ which completes the proof.
\end{proof}

Let
$$\begin{array}{ll}
G_t&:=P(\ttau>t|\scr
F_t)=\dint_t^\infty\dint_\bb Rp_t(s,l)\eta(ds,dl)\quad\text{and}\\
G^*_t&:=P^*(\ttau>t|\scr
F_t)=P^*(\ttau>t)=P(\ttau>t)=\dint_0^t\dint_\bb R\eta(ds,dl),
\end{array}$$
from
Callegaro, Jeanblanc and Zargari(2010)$^{\text{\cite{Callegaro-Jeanblanc-Zargari(2010)}}}$,
we find that $G$ is a $(P,\bb F)$-supermartingale and $G^*$ is a
deterministic continuous and decreasing function.

\begin{theorem}\label{theorem-Callegaro et al(2010)}
Let $y_t(\ttau,L)$ be a $\scr G^{\ttau,L}_t$-measurable
$P$-integrable r.v., then for $s\le t$,
$$E(y_t(\ttau,L)|\scr G_s)=\wt y_s\bb I_{s<\ttau}+\wh y_s(\ttau,L)\bb I_{\ttau\le s}$$
with
$$\begin{array}{ll}
\wt y_s&=\dfrac{1}{G_s}E\Big(\dint_s^\infty\int_\bb
Ry_t(u,l)p_t(u,l)\eta(du,dl)\Big|\scr
F_s\Big)\\
\wh y_s(u,l)&=\dfrac{1}{p_s(u,l)}E\big\{y_t(u,l)p_t(u,l)\big|\scr
F_s\big\}.
\end{array}$$
\end{theorem}
\begin{proof}
The proof follows from the proof of Lemma 1.5 of
Callegaro-Jeanblanc-Zargari(2010)(\cite{Callegaro-Jeanblanc-Zargari(2010)}).
\end{proof}

\subsection{The $(P,\bb F)$-density process}
Since for any $s\geq 0,l\in \bb R$, $p(s,l)=\{p_t(s,l);t\ge 0\}$ is a strictly
positive $(P,\bb F)$-martingale which implies that $p_t(s,l)$ can be
represented in the following form
\begin{equation}\label{eqn-p-2}
\begin{array}{ll}
p_t(s,l)=E(p_s(s,l)|\scr F_t)\bb
I_{t<s}+p_s(s,l)\exp\bigg\{\dint_s^t\theta_1(u;s,l)'d
W(u)-\dfrac{1}{2}\dint_s^t\|\theta_1(u;s,l)\|^2du\bigg\}\\
\qquad\times
\exp\bigg\{\dint_s^t\dint_E\big\{\ln(1+\theta_2(u,y;s,l))\big\}\mu(du,dy)
-\dint_s^t\dint_E\theta_2(u,y;s,l)\nu(du,dy)\bigg\}\bb I_{t\ge
s},\end{array}
\end{equation}
where $\theta_{1}(\cdot;s,l)$ is an $\bb R^d$-valued $\bb F$-predictable process with $\theta_{1}(\cdot;s,l)\in L^2_{loc}(W)$ and $\theta_{2}(\cdot,\cdot;s, l)$
 is a $\widetilde{\scr P}(\bb F)$-measurable function with $\theta_{2}(\cdot,\cdot;s, l)\in G_{loc}(\mu) $.
Therefore, the conditional density $p_t(s,l)$ is completely determined by
$p_s(s,l)$, $\theta_1(u;s,l)\bb I_{u>s}$ and $\theta_2(u,y;s,l)\bb
I_{u>s}$.
For given $(p_s(s,l), \theta_1, \theta_2)$,  we define
\begin{equation} \label{Z}
\begin{array}{ll}
Z^{\theta_1,\theta_2}_{s,l}(t):=\exp\bigg\{\dint_s^t\theta_1(u;s,l)'d
W(u)-\dfrac{1}{2}\dint_s^t\|\theta_1(u;s,l)\|^2du\bigg\}\\
\qquad\times\exp\bigg\{\dint_s^t\dint_E\big\{\ln(1+\theta_2(u,y;s,l))\big\}\mu(du,dy)
-\dint_s^t\dint_E\theta_2(u,y;s,l)\nu(du,dy)\bigg\},\end{array}
\end{equation}
then the conditional density $p_t(s,l)$ has following representation:
\begin{equation} \label{Z}
\begin{array}{ll}
p_t(s,l)=E(p_s(s,l)|\scr F_t)\bb
I_{t<s}+p_s(s,l)Z^{\theta_1,\theta_2}_{s,l}(t)\bb I_{t\ge
s}.\end{array}
\end{equation}

For parameters of $p_s(s,l)$, $\theta_1$, $\theta_2$,
we have the following theorem.

\begin{theorem}\label{theorem-p(s;s,l)}
 If $p_t(s,l)$
 is the density process of a pair $(\ttau,L)$ with
respect to $(P,\bb F)$ and has the representation as (\ref{Z}), then $p_s(s,l)$, $\theta_1(u;s,l)\bb
I_{u>s}$ and $\theta_2(u,y;s,l)\bb I_{u>s}$  satisfy the following
equation
\begin{equation}\label{condition-p(s;s,l)}
\begin{array}{ll}\dint_0^\infty\dint_\bb Rp_s(s,l)\bb
\eta(ds,dl)=1-\dint_0^\infty\Big\{\dint_0^u\dint_\bb
Rp_s(s,l)Z^{\theta_1,\theta_2}_{s,l}(u-)\theta_1(u;s,l)\eta(ds,dl)\Big\}'d
W(u)\\\qquad-\dint_0^\infty\dint_E\Big\{\dint_0^u\dint_\bb
Rp_s(s,l)Z^{\theta_1,\theta_2}_{s,l}(u-) \theta_2(u,y;s,l))
\eta(ds,dl)\Big\}\{\mu(du,dy)-\nu(du,dy)\}.\end{array}
\end{equation}
\end{theorem}

\begin{proof}
We first note that the following equation by Fubini theorem.
\begin{equation}\label{eqn-Fubini-1}
\dint_t^\infty\dint_\bb RE(p_s(s,l)|\scr F_t)\eta(ds,dl)
=E\Big(\dint_t^\infty\dint_\bb Rp_s(s,l)\bb \eta(ds,dl)\Big|\scr
F_t\Big).
\end{equation}
When $t=0$, the above equation becomes
\begin{equation}
E\Big(\dint_0^\infty\dint_\bb Rp_s(s,l)\bb \eta(ds,dl)
\Big)=\dint_0^\infty\dint_\bb Rp_0(s,l)\eta(ds,dl)=1.
\end{equation}

Since
$$\begin{array}{ll}
Z^{\theta_1,\theta_2}_{s,l}(t)=1+\dint_s^tZ^{\theta_1,\theta_2}_{s,l}(u-)\theta_1(u;s,l)'d
W(u)\\\qquad\qquad\qquad+\dint_s^t\dint_EZ^{\theta_1,\theta_2}_{s,l}(u-)\theta_2(u,y;s,l))\{\mu(du,dy)-\nu(du,dy)\},\end{array}$$
then from $\dint_0^\infty\int_\bb Rp_t(s,l)\eta(ds,dl)\equiv1$ for each
$t$, one can see that

$$\begin{array}{ll}
1=\dint_0^\infty\dint_\bb R p_t(s,l)\eta(ds,dl)\\
=\dint_t^\infty\dint_\bb R p_t(s,l)\eta(ds,dl)+\dint_0^t\dint_\bb R p_t(s,l)\eta(ds,dl)\\
=E\Big(\dint_t^\infty\dint_\bb R p_s(s,l)\eta(ds,dl)\Big|\scr
F_t\Big)+\dint_0^t\dint_\bb
Rp_s(s,l)Z^{\theta_1,\theta_2}_{s,l}(t)
\eta(ds,dl)\\
=E\Big(\dint_0^\infty\dint_\bb Rp_s(s,l)\bb \eta(ds,dl)\Big|\scr
F_t\Big)\\
\qquad-\dint_0^t\dint_\bb Rp_s(s,l)\bb
\eta(ds,dl)+\dint_0^t\dint_\bb
Rp_s(s,l)Z^{\theta_1,\theta_2}_{s,l}(t)
\eta(ds,dl)\\

=E\Big(\dint_0^\infty\dint_\bb Rp_s(s,l)\bb \eta(ds,dl)\Big|\scr
F_t\Big)\\
\qquad+\dint_0^t\dint_\bb
Rp_s(s,l)\Big\{\dint_s^tZ^{\theta_1,\theta_2}_{s,l}(u-)\theta_1(u;s,l)'d
W(u)\\\qquad\qquad\qquad
+\dint_s^t\dint_EZ^{\theta_1,\theta_2}_{s,l}(u-)
\theta_2(u,y;s,l))\{\mu(du,dy)-\nu(du,dy)\}\Big\}
\eta(ds,dl)\\

=E\Big(\dint_0^\infty\dint_\bb Rp_s(s,l)\bb \eta(ds,dl)\Big|\scr
F_t\Big)\\
\qquad+\dint_0^t\Big\{\dint_0^u\dint_\bb
Rp_s(s,l)Z^{\theta_1,\theta_2}_{s,l}(u-)\theta_1(u;s,l)\eta(ds,dl)\Big\}'d
W(u)\\\qquad+\dint_0^t\dint_E\Big\{\dint_0^u\dint_\bb
Rp_s(s,l)Z^{\theta_1,\theta_2}_{s,l}(u-) \theta_2(u,y;s,l))
\eta(ds,dl)\Big\}\{\mu(du,dy)-\nu(du,dy)\},
\end{array}$$
where the last equality we have used the stochastic Funini theorem for general semimartingales,
see Theorem 4.1.1 of Jeanblanc, Yor and Chesney(2009)
$^{\text{\cite{Jeanblanc-Yor-Chesney2009}}}$.

 Since $E\Big(\dint_0^\infty\dint_\bb Rp_s(s,l)\bb \eta(ds,dl)\Big|\scr
F_t\Big)$ is an uniformly integrable $(P,\bb F)$-martingales, taking $t\longrightarrow\infty$ and
hence equation (\ref{condition-p(s;s,l)}) holds.
\end{proof}

More interesting, we can consider the inverse problem, given a law $\eta(s,l)$ with
 \begin{equation}\label{eta}
\begin{array}{ll}
 \dint_0^\infty\dint_\bb R\eta(ds,dl)=1
 \end{array}
\end{equation}
and $(p_s(s,l),\theta_1,\theta_2)$,
whether can we construct a ``conditional density"?  The following theorem gives the answer.

\begin{theorem}\label{sufficient conditions-p_t-1}
 Given $(p_s(s,l),\theta_1,\theta_2)$ satisfying equation (\ref{condition-p(s;s,l)}) such that $E(\dint_0^\infty\int_\bb R p_s(s,l)\eta(ds,dl))=1$ and
 $Z^{\theta_1,\theta_2}_{s,l}(t)$ is a $(P, \bb F)$ martingale,
 then there exists a pair $(\ttau^*,L^*)$, such that $\widetilde{\tau}^*$ is a non-negative random variable and $L^*$ is a
random variable on $(\Omega,\scr F)$, and
$$\begin{array}{ll}
p_t(s,l):=E(p_s(s,l)|\scr F_t)\bb
I_{t<s}+p_s(s,l)Z^{\theta_1,\theta_2}_{s,l}(t)\bb I_{t\ge
s}\end{array}$$ is the density process of  $(\ttau^*,L^*)$.
\end{theorem}

\begin{proof}
First, it is noted that $p_t(s, l)$ is an $(P, \bb F)$ martingale for any
$s\geq 0, l\in \bb R$ and
from the proved process of Theorem \ref{theorem-p(s;s,l)}, one can see that for any
$t\geq 0$,
$$\begin{array}{ll}1=\dint_t^\infty\dint_\bb R E(p_s(s,l)|\scr F_t)\eta(ds,dl)+ \dint_0^t\dint_\bb Rp_s(s,l)Z^{\theta_1,\theta_2}_{s,l}(t)\eta(ds,dl)\\
=\dint_0^\infty\dint_\bb Rp_t(s,l)\eta(ds,dl).
\end{array}$$
Particularly,  from equation
 (\ref{eta}), one has $p_0(s, l)=1$, for any $t\geq 0, l\in \bb R$.

Since $ \dint_0^\infty\dint_\bb R\eta(ds,dl)=1$, then we can construct a probability space
 $(\Omega', P')$ and
a non-negative random
variable $\ttau^*$ and  a
random variable $L^*$ on $\bb F$ such that
the law of $(\ttau^*, L^*)$ is
$$P'((\ttau^*,L^*)\in B)=\int_B\eta(ds,dl)\quad \text{for any}\quad  B\in\scr B(\bb R^+\times\bb R).$$
Define the product space
$$\Omega^*=\Omega\times\Omega', ~\scr G^*=\scr F_\infty\otimes\sigma(\ttau, L)$$
and the filtration
$$\scr G_t^*=\scr F_t\otimes\sigma(\ttau, L).$$
Define the measure $\overline{P}$ on $(\Omega^*, \scr G^*)$
$$\overline{P}=P\mid_{\scr F_t}\otimes P'\mid_{\sigma(\ttau, L)},$$
then
$(\ttau^*,L^*)$ is independent of $\scr F_t$  on $(\Omega^*, \scr G^*, P)$
and
$\overline{P}|_{\scr F_t}=P|_{\scr F_t}$
and $\overline{P}|_{\sigma(\ttau,L)}=P|_{\sigma(\ttau,L)}$,
One also checks that $p_t(\ttau^*, L^*)$ is an $\scr G_t^*$ martingale
(similar to
Grorud and Pontier(1998)$^{\text{\cite{Grorud-Pontier-1998}}}$  or
Amendinger(1999)$^{\text{\cite{Amendinger(1999)}}}$).
Thus one can define a locally
equivalent probability measure $Q$ by
 $\dfrac{dQ}{d\overline{P}}|_{\scr G_t^*}=p_t(\ttau^*, L^*)$,
then
for any $B\in\scr B(\bb R^+\times\bb R)$ and for every $t\ge 0$,  from Bayes' formula,  one obtains
$$
\begin{array}{ll}
Q((\ttau^*,L^*)\in B|\scr F_t)\\
=\dfrac{E(p_t(\ttau^*,L^*)\bb I_{(\ttau^*,L^*)\in B}|\scr F_t)}{E(p_t(\ttau^*,L^*)|\scr F_t)}\\
=\int_Bp_t(s,l)\eta(ds,dl).
\end{array}
$$
then $p_t(s, l)$ is the conditional
density of $(\ttau^*,L^*)$ with respect to $\eta$.
\end{proof}

\begin{theorem}\label{sufficient conditions-p_t-2} If we
 given a positive
 $\scr O(\bb F)\otimes\scr B(\bb R)$-measurable function
 $p_s(s,l)$ satisfying $E(\dint_0^\infty\int_\bb R p_s(s,l)\eta(ds,dl))=1$ and
 $E(\dint_t^\infty\int_\bb R p_s(s,l)\eta(ds,dl)\mid\scr F_t)<1$ for any $t>0$,
  then we can find
 a pair $(\theta_1, \theta_2)$ such that
$(p_s(s,l),\theta_1,\theta_2)$ satisfying equation (\ref{condition-p(s;s,l)}).
\end{theorem}

\begin{proof}
For any $t\geq 0$,
let $Y_t:=1-E(\dint_t^\infty\int_\bb R p_s(s,l)\eta(ds,dl)\mid\scr F_t)$, then
$Y_t$ is a $(P, \bb F)$ positive submartingale and
$Y_t-\dint_0^t\int_\bb R p_s(s,l)\eta(ds,dl)
=1-E(\dint_0^\infty\int_\bb R p_s(s,l)\eta(ds,dl)\mid\scr F_t)$
is a  $(P, \bb F)$ martingale.
 By
Assumption \ref{assum-filtration-F}
, there exists a pair
$(\widetilde{\theta}_1, \widetilde{\theta}_2)$ such that

 \begin{equation}\label{Y-01}
\begin{array}{ll}
Y_t=\dint_0^t\int_\bb R p_s(s,l)\eta(ds,dl)-\dint_0^tY_{u-}\widetilde{\theta}_1(u)'dW(u)-\dint_0^t\dint_EY_{u-}\widetilde{\theta}_2(u, y)
\{\mu(du,dy)-\nu(du,dy)\},
 \end{array}
\end{equation}
where $\wt \theta_{1}(\cdot)$ is an $\bb R^d$-valued $\bb F$-predictable process with $\wt \theta_{1}(\cdot)\in L^2_{loc}(W)$ and $\wt\theta_{2}(\cdot,\cdot)$
 is a $\widetilde{\scr P}(\bb F)$-measurable function with $\wt\theta_{2}(\cdot,\cdot)\in G_{loc}(\mu)$.
 For any
$s\geq 0, l\in \bb R$, we define $\theta_1(u;s,l)$
and $\theta_2(u,y;s,l)$
$$\left\{\begin{array}{ll}
\theta_1(u;s,l)&:=-\wt\theta_1(u)\bb I_{u>s}\\
\theta_2(u,y;s,l)&:=-\wt\theta_2(u,y)\bb I_{u>s}.
\end{array}\right.$$
One can see from the defintion of  $Z^{\theta_1,\theta_2}_{s,l}$,  $Z^{\theta_1,\theta_2}_{s,l}$ is the solution of following BSDE
$$\begin{array}{ll}
Z^{\theta_1,\theta_2}_{s,l}(t)
=1-\dint_s^tZ^{\theta_1,\theta_2}_{s,l}(u-)\wt\theta_1(u)'d
W(u)\\\qquad\qquad\qquad-\dint_s^t\dint_EZ^{\theta_1,\theta_2}_{s,l}(u-)
\wt\theta_2(u,y)\{\mu(du,dy)-\nu(du,dy)\}.
\end{array}$$
Let $\wt Y_t:=\dint_0^t\dint_\bb R p_s(s, l)Z^{\theta_1,\theta_2}_{s,l}(t)$ then, $\wt Y_t$ satisfies
$$\begin{array}{ll}
\wt Y_t=\dint_0^t\dint_\bb R p_s(s, l)Z^{\theta_1,\theta_2}_{s,l}(t)\eta(ds,dl)\\
=\dint_0^t\dint_\bb R p_s(s, l)\eta(ds,dl)-\dint_0^t\dint_\bb R p_s(s, l)\dint_s^tZ^{\theta_1,\theta_2}_{s,l}(u-)\wt\theta_1(u)'d
W(u)\eta(ds,dl)\\\qquad\qquad\qquad-\dint_0^t\dint_\bb R p_s(s, l)\dint_s^t\dint_EZ^{\theta_1,\theta_2}_{s,l}(u-)
\wt\theta_2(u,y)\{\mu(du,dy)-\nu(du,dy)\}\eta(ds,dl)\\
=\dint_0^t\dint_\bb R p_s(s, l)\eta(ds,dl)-\dint_0^t\Big\{\dint_0^u\dint_\bb R p_s(s, l)Z^{\theta_1,\theta_2}_{s,l}(u-)\eta(ds,dl)\Big\}\wt\theta_1(u)'d
W(u)\\\qquad\qquad\qquad-\dint_0^t\dint_E\Big\{\dint_0^u\dint_\bb R p_s(s, l)Z^{\theta_1,\theta_2}_{s,l}(u-)
\eta(ds,dl)\Big\}\wt\theta_2(u,y)\{\mu(du,dy)-\nu(du,dy)\}\\
=\dint_0^t\int_\bb R p_s(s,l)\eta(ds,dl)-\dint_0^t\wt Y_{u-}\widetilde{\theta}_1(u)'dW(u)-\dint_0^t\dint_E\wt Y_{u-}\widetilde{\theta}_2(u, y)
\{\mu(du,dy)-\nu(du,dy)\},
\end{array}$$
From the unique solution of SDE, we obtain
$Y_t=\wt Y_t
=\dint_0^t\dint_\bb R p_s(s,l)Z^{\theta_1,\theta_2}_{s,l}(t)\eta(ds,dl), a.s.$.
From
equation (\ref{Y-01}) and let $t=\infty$, we know that
\begin{equation}
\begin{array}{ll}\dint_0^\infty\dint_\bb Rp_s(s,l)\bb
\eta(ds,dl)=1-\dint_0^\infty\Big\{\dint_0^u\dint_\bb
Rp_s(s,l)Z^{\theta_1,\theta_2}_{s,l}(u-)\theta_1(u;s,l)\eta(ds,dl)\Big\}'d
W(u)\\\qquad-\dint_0^\infty\dint_E\Big\{\dint_0^u\dint_\bb
Rp_s(s,l)Z^{\theta_1,\theta_2}_{s,l}(u-) \theta_2(u,y;s,l))
\eta(ds,dl)\Big\}\{\mu(du,dy)-\nu(du,dy)\},\end{array}
\end{equation}
which implies that $(p_s(s,l),\theta_1,\theta_2)$ satisfying equation (\ref{condition-p(s;s,l)})
\end{proof}

The following corollary gives an explicit expression of the survival process $G_t$ by $(p_s(s,l),\theta_1,\theta_2)$.

\begin{corollary}
Under the conditions of Theorem \ref{theorem-p(s;s,l)}, the survival
process of $\ttau$ with respect to $(P,\bb F)$ is given by
\begin{equation}\tag{\ref{eqn-G-Doob-Meyer}}
\begin{array}{ll}G_t
&=1-\dint_0^t\int_\bb Rp_s(s,l)\eta(ds,dl)\\
&\qquad-\dint_0^t\Big\{\dint_0^u\dint_\bb
Rp_s(s,l)Z^{\theta_1,\theta_2}_{s,l}(u-)\theta_1(u;s,l)\eta(ds,dl)\Big\}'d
W(u)\\
&\qquad-\dint_0^t\dint_E\Big\{\dint_0^u\dint_\bb
Rp_s(s,l)Z^{\theta_1,\theta_2}_{s,l}(u-) \theta_2(u,y;s,l))
\eta(ds,dl)\Big\}\{\mu(du,dy)-\nu(du,dy)\}.
\end{array}
\end{equation}
\end{corollary}

\begin{proof}
From the definition of $G_t$, one sees that
$$\begin{array}{ll}G_t
&=P(\ttau>t|\scr F_t)=\dint_t^\infty\int_\bb Rp_t(s,l)\eta(ds,dl)\\
&=\dint_t^\infty\int_\bb RE[p_s(s,l)|\scr F_t]\eta(ds,dl)\\
&=E\Big[\dint_t^\infty\int_\bb Rp_s(s,l)\eta(ds,dl)\Big|\scr F_t\Big]\\
&=-\dint_0^t\int_\bb
Rp_s(s,l)\eta(ds,dl)+E\Big[\dint_0^\infty\int_\bb
Rp_s(s,l)\eta(ds,dl)\Big|\scr F_t\Big]\\
&=1-\dint_0^t\int_\bb Rp_s(s,l)\eta(ds,dl)\\
&\qquad-\dint_0^t\Big\{\dint_0^u\dint_\bb
Rp_s(s,l)Z^{\theta_1,\theta_2}_{s,l}(u-)\theta_1(u;s,l)\eta(ds,dl)\Big\}'d
W(u)\\
&\qquad-\dint_0^t\dint_E\Big\{\dint_0^u\dint_\bb
Rp_s(s,l)Z^{\theta_1,\theta_2}_{s,l}(u-) \theta_2(u,y;s,l))
\eta(ds,dl)\Big\}\{\mu(du,dy)-\nu(du,dy)\},
\end{array}$$
which completes the proof.
\end{proof}

We give a special example as the end of this subsection.
\begin{example}\label{example}
If
$$\dint_0^\infty\dint_\bb Rp_s(s,l)\eta(dl)ds=1,\quad a.s.,$$
let $\theta_1(u;s,l)=0$ and $\theta_2(u,y;s,l):=0$, then one can see that $p_s(s,l)$, $\theta_1(u;s,l)$, $\theta_2(u;s,l)$ satisfies (\ref{condition-p(s;s,l)}). Hence $
p_t(s,l):=E(p_s(s,l)|\scr F_t)\bb
I_{t<s}+p_s(s,l)Z^{\theta_1,\theta_2}_{s,l}(t)\bb I_{t\ge
s}$ is the density process of some $(\ttau,L)$. In this case,
$$\begin{array}{ll}
G_t&=1-\dint_0^t\int_\bb Rp_s(s,l)\eta(dl)ds.
\end{array}$$
\end{example}

\subsection{$\bb G$-martingales' characterization}

Similar to the proof of
Proposition 2.2 of Callegaro, Jeanblanc and Zargari(2010)$^{\text{
\cite{Callegaro-Jeanblanc-Zargari(2010)}}}$, one can show the following theorem

\begin{theorem}
[Characterization of $(P,\bb G)$-martingales in terms of $(P, \bb
F)$-martingales]\label{theorem-G-martingale character} Let
$y=\{y_t;t\ge0\}$, where $y_t:= \wt y_t\bb I_{t<\ttau} +\wh
y_t(\ttau,L)\bb I_{t\ge \ttau}$, be a $\bb G$-adapted process, then $y$
is a  $(P,\bb G)$-martingale if and only
if the following two conditions are satisfied:
\begin{itemize}
  \item[(i)] for $\eta$-almost every $u\ge 0$ and $l\in \bb R$, $\{\wh y_t(u; l)p_t(u; l); t\ge u\}$ is a
  $(P, \bb F)$-martingale;
  \item[(ii)]the process $\{\wt y_tG_t + \int_0^t\int_\bb R \wh y_u(u, l)p_u(u, l)\eta(du,dl); t \ge
  0\}$
is a $(P,\bb F)$-martingale.
\end{itemize}
\end{theorem}

Similarly, we have the following corollary
\begin{corollary}
\label{theorem-(P*,G)-martingale character} Let $y=\{y_t;t\ge0\}$,
where $y_t:= \wt y_t\bb I_{t<\ttau} +\wh y_t(\ttau,L)\bb I_{t\ge \ttau}$, be
a $\bb G$-adapted process, then $y$ is a
$(P^*,\bb G)$-martingale if and only if the following two conditions
are satisfied:
\begin{itemize}
  \item[(i)] for $\eta$-almost every $u\ge 0$ and $l\in \bb R$, $\{\wh y_t(u; l); t\ge u\}$ is a
  $(P^*, \bb F)$-martingale;
  \item[(ii)]the process $\{\wt y_tG^*_t + \int_0^t\int_\bb R \wh y_u(u, l)\eta(du,dl); t \ge
  0\}$
is a $(P^*,\bb F)$-martingale.
\end{itemize}
\end{corollary}

\section{Canonical decomposition of a $(P,\bb F)$-martingale in $(P,\bb G)$}
We now consider  the canonical decomposition of any $(P,\bb F)$
martingale $m$ in the enlarged filtration $\bb G$ respectively under
Assumption \ref{assumption-01}. Form Theorem \ref{theorem-p(s;s,l)}, one
can see that $p_t(s,l)$ can be determined by $p_s(s,l)$,
$\theta_1(u;s,l)\bb I_{u>s}$ and $\theta_2(u,y;s,l)\bb I_{u>s}$. We have
$$p_t(s,l)=E[p_s(s,l)|\scr F_t]\bb
I_{t<s}+p_s(s,l)Z^{\theta_1,\theta_2}_{s,l}(t)\bb I_{t\ge s}$$ and
$$\begin{array}{ll}\\
p_t(s,l)\bb I_{t\ge s}=p_s(s,l)\bb I_{t\ge s}+\dint_0^tp_{u-}(s,l)\theta_1(u;s,l)'\bb I_{u>s}dW(u)\\
\qquad +\dint_0^t\int_Ep_{u-}(s,l)\theta_2(u,y;s,l)\bb
I_{u>s}\{\mu(du,dy)-\nu(du,dy)\},
\end{array}$$
thus
$$\begin{array}{ll}
p_t(\ttau,L)\bb I_{t\ge \ttau}=p_\ttau(\ttau,L)\bb I_{t\ge \ttau}+\dint_0^tp_{u-}(\ttau,L)\theta_1(u;\ttau,L)'\bb
I_{u>\ttau}dW(u)\\
\qquad +\dint_0^t\int_Ep_{u-}(\ttau,L)\theta_2(u,y;\ttau,L)\bb
I_{u>\ttau}\{\mu(du,dy)-\nu(du,dy)\}.
\end{array}$$\\
Recall that $G^*_t=P(\ttau>t)$ is a deterministic continuous
function satisfying $0<G^*_{t}<1$ for each $t\in(0,\infty)$, since there
are no atoms.
To obtain  the canonical decomposition of a $(P,\bb F)$
martingale in the filtration $\bb G$, we need the following lemma:

\begin{lemma}\label{lemma-Af}
For any positive $\scr O(\bb F)\times\scr B$-measurable function
$f(s,l)$ such that $E[\dint_{0}^{\infty}\int_\bb Rf(s,l)$
$\eta(ds,dl)]<\infty$, let
$$A^{f,*}_t:=\int_0^t\int_\bb R \dfrac{f(s,l)}{G^*_{s-}}\eta(ds,dl),$$
then $A^{f,*}$ is a continuous increasing $\bb F$-adapted process
and
$$M^{f,*}_t=f(\ttau,L)\bb I_{t\ge\ttau}-A^{f,*}_{t\wedge\ttau}$$
is a $(P^*,\bb G)$-martingale, i.e., $A^{f,*}_{\0\wedge\ttau}$ is
the $(P^*,\bb G)$-compensator of $f(\ttau,L)\bb I_{t\ge\ttau}$.
\end{lemma}
\begin{proof}
For any $t_1<t_2$, we have
\begin{equation}\label{eqn-independent}
\begin{array}{ll}
E_{P^*}\big[f(\ttau,L)\bb I_{\ttau\le t}\big|\scr F_{t}\big]
=\dint_{0}^{t}\int_\bb Rf(s,l)\eta(ds,dl).
\end{array}
\end{equation}
In fact, one can see from the independence of $(\ttau,L)$ with
respect to $\scr F_t$ under $P^*$ that (\ref{eqn-independent}) holds
for all positive $\scr B(\bb R^+\times\bb R)$-measurable functions
$f(s,l)$. In general, one gets from the independence lemma (see
Lemma 2.3.4 of Shreve(2003)$^{\text{\cite{Shreve2003}}}$) and the
monotone class theorem that (\ref{eqn-independent}) still holds for
any positive $\scr O(\bb F)\times\scr B$-measurable function
$f(s,l)$. And since there are no atoms, one can see that
$\dint_{0}^{t}\int_\bb Rf(s,l)\eta(ds,dl)$ is a continuous
increasing $\bb F$-adapted process, thus a $\bb F$-predictable
process.\par

For any $t_1<t_2$, we have
$$\begin{array}{ll}
E_{P^*}\big[f(\ttau,L)\bb I_{t_1<\ttau\le t_2}\big|\scr F_{t_1}\big]
&=E_{P^*}\big[\dint_{t_1}^{t_2}\int_\bb Rf(s,l)\eta(ds,dl)\big|\scr
F_{t_1}\big]
\end{array}$$
and
$$\begin{array}{ll}
E_{P^*}\Big[\dint_{t_1}^{t_2}\int_\bb R\dfrac{\bb
I_{s\le\ttau}}{G^*_s}f(s,l)
\eta(ds,dl)\Big|\scr F_{t_1}\Big]\\
\qquad=E_{P^*}\Big[\dint_{t_1}^{t_2}\int_\bb R\dfrac{E_{P^*}[\bb
I_{s\le\ttau}|\scr F_s]}{G^*_{s-}}f(s,l)
\eta(ds,dl)\Big|\scr F_{t_1}\Big]\\
\qquad=E_{P^*}\Big[\dint_{t_1}^{t_2}\int_\bb
Rf(s,l)\eta(ds,dl)\Big|\scr F_{t_1}\Big],
\end{array}$$
hence
$$\begin{array}{ll}
E_{P^*}[M^{f,*}_{t_2}-M^{f,*}_{t_1}|\scr
G_{t_1}]=E_{P^*}\Big[f(\ttau,L)\bb I_{t_1<\ttau\le
t_2}-\dint_{t_1}^{t_2}\int_\bb
R\dfrac{\bb I_{s\le\ttau}}{G^*_s}f(s,l)\eta(ds,dl)\Big|\scr G_{t_1}\Big]\\
=\dfrac{1}{G^*_{t_1}}\Big\{E_{P^*}\big[f(\ttau,L)\bb I_{t_1<\ttau\le
t_2}\big|\scr F_{t_1}\big]\\
\qquad\qquad-E_{P^*}\Big[\dint_{t_1}^{t_2}\int_\bb R\dfrac{\bb
I_{s\le\ttau}}{G^*_s}f(s,l)
\eta(ds,dl)\Big|\scr F_{t_1}\Big]\Big\}\bb I_{t_1<\ttau}\\
=0,
\end{array}$$
thus $E_{P^*}[M^{f,*}_{t_2}|\scr G_{t_1}]=M^{f,*}_{t_1}$ and
$M^{f,*}$ is a $(P^*,\bb G)$-martingale.
\end{proof}

It is noted that
$\dfrac{1}{p_t(\ttau,L)}$ is the density process of $P^*$ with
respect to $(P,\bb G^{\ttau,L})$, one can see that the density
process of $P^*$ with respect to $(P,\bb G)$ is given by
$$\begin{array}{ll}
L^*_t:=E\Big[\dfrac{1}{p_t(\ttau,L)}\Big|\scr G_t\Big]\\
\qquad=\dfrac{1}{G_t}E\Big(\dint_t^\infty\int_\bb
R\dfrac{1}{p_t(u,l)}p_t(u,l)\eta(du,dl)\Big|\scr F_s\Big)\bb
I_{t<\ttau}+\dfrac{1}{p_t(\ttau,L)}\bb I_{t\ge\ttau}\\
\qquad=\dfrac{1}{G_t}E\Big(\dint_t^\infty\int_\bb
R\eta(du,dl)\Big|\scr F_s\Big)\bb
I_{t<\ttau}+\dfrac{1}{p_t(\ttau,L)}\bb I_{t\ge\ttau}\\
\qquad=\dfrac{G^*_t}{G_t}\bb I_{t<\ttau}+\dfrac{1}{p_t(\ttau,L)}\bb
I_{t\ge\ttau},
\end{array}$$
thus
$\dfrac{1}{L^*_t}
=\dfrac{G_t}{G^*_t}\bb I_{t<\ttau}+p_t(\ttau,L)\bb
I_{t\ge\ttau}$
 is a  $(P^*,\bb
G)$-martingale, then we
have:

\begin{corollary}
Let
$$N_1(t):=p_\ttau(\ttau,L)\bb I_{t\ge \ttau}-\dint_0^{t\wedge\ttau}\int_\bb
R\dfrac{1}{G^*_s}p_s(s,l)\eta(ds,dl),$$
$$N_2(t):=\dfrac{G_{\ttau}}{G^*_\ttau}\bb
I_{t\ge\ttau}+\dint_0^{t\wedge\ttau}\dfrac{G_{u}}{(G^*_{u})^2}dG^*_u,$$
 then $\{N_1(t);t\ge 0\}$ and $\{N_2(t);t\ge 0\}$ are
 $(P^*,\bb G)$-martingales.
\end{corollary}

\begin{proof}
 One can see that
$$\begin{array}{ll}N_2(t)
&=\dfrac{G_{\ttau}}{G^*_\ttau}\bb
I_{t\ge\ttau}+\dint_0^{t\wedge\ttau}\dfrac{G_{s}}{(G^*_{s})^2}dG^*_s\\
&=\dfrac{G_{\ttau}}{G^*_\ttau}\bb
I_{t\ge\ttau}-\dint_0^{t\wedge\ttau}\dint_\bb
R\dfrac{1}{G^*_s}\dfrac{G_{s}}{G^*_{s}}\eta(ds,dl),
\end{array}$$
since $G^*_s=P^*(\ttau>s)=\dint_s^\infty\int_\bb R\eta(du,dl)$.
It is noted that $E(G_{s})=E(\dint_s^\infty\dint_\bb Rp_s(u,l)\eta(du,dl))=\dint_s^\infty\dint_\bb RE(p_s(u,l))\eta(du,dl)
=\dint_s^\infty\dint_\bb R\eta(du,dl)=G^*_s$ and
$$\begin{array}{ll}
E(\dint_{0}^{\infty}\int_\bb R\dfrac{G_{s}}{G^*_s}
\eta(ds,dl))\\
=\dint_{0}^{\infty}\int_\bb R\dfrac{1}{G^*_s}E(G_{s})
\eta(ds,dl)\\
=\dint_{0}^{\infty}\int_\bb R
\eta(ds,dl)\\
=1<\infty
\end{array}$$ and
$E(p_\ttau(\ttau, L))=E(E(p_\ttau(\ttau, L)|\scr F_t))=E(\dint_{0}^{\infty}\int_\bb Rp_s(s, l)
\eta(ds,dl))=1<\infty$.
Then from
Lemma \ref{lemma-Af}, one can see that $\{N_1(t);t\ge 0\}$ and $\{N_2(t);t\ge 0\}$ are
 $(P^*,\bb G)$-martingales.
\end{proof}

\begin{theorem}\label{prop-presentation of Z*}
Let $Z^*_t:=\dfrac{G_t}{G^*_t}\bb I_{t<\ttau}+p_t(\ttau,L)\bb
I_{t\ge\ttau}$, then $Z^*$ is a  $(P^*,\bb
G)$-martingale with the following decomposition
$$\begin{array}{ll}Z^*_t
&=1-\dint_0^tZ^*_{u-}\bigg\{\dfrac{1}{G_{u-}}\dint_0^u\dint_\bb
Rp_s(s,l)Z^{\theta_1,\theta_2}_{s,l}(u-)\theta_1(u;s,l)\eta(ds,dl)\bb
I_{u\le\ttau}\\
&\qquad\qquad\qquad-\theta_1(u;\ttau,L)\bb I_{u>\ttau}\bigg\}'dW(u)\\
&\qquad-\dint_0^{t}\dint_EZ^*_{u-}\bigg\{\dfrac{1}{G_{u-}}\dint_0^u\dint_\bb
Rp_s(s,l)Z^{\theta_1,\theta_2}_{s,l}(u-) \theta_2(u,y;s,l))
\eta(ds,dl)\bb I_{u\le\ttau}\\
&\qquad\qquad\qquad -\theta_2(u,y;\ttau,L)\bb
I_{u>\ttau}\bigg\}\{\mu(du,dy)-\nu(du,dy)\}\\
&\qquad+N_1(t)-N_2(t).
\end{array}$$
\end{theorem}
\begin{proof}
One can see that
$Z^*_t=\dfrac{1}{p_t(\ttau,L)}=\dfrac{1}{E\Big[\dfrac{1}{p_t(\ttau,L)}\Big|\scr G_t\Big]}$,
from which $Z^*$ is a
$(P^*,\bb G)$-martingale. Since
$$\begin{array}{ll}\dfrac{G_t}{G^*_t}
&=1-\dint_0^t\int_\bb R\dfrac{1}{G^*_s}p_s(s,l)\eta(ds,dl)\\
&\qquad-\dint_0^t\dfrac{1}{G^*_u}\Big\{\dint_0^u\dint_\bb
Rp_s(s,l)Z^{\theta_1,\theta_2}_{s,l}(u-)\theta_1(u;s,l)\eta(ds,dl)\Big\}'d
W(u)\\
&\qquad-\dint_0^t\dint_E\dfrac{1}{G^*_u}\Big\{\dint_0^u\dint_\bb
Rp_s(s,l)Z^{\theta_1,\theta_2}_{s,l}(u-) \theta_2(u,y;s,l))
\eta(ds,dl)\Big\}\{\mu(du,dy)-\nu(du,dy)\}\\
&\qquad-\dint_0^t\dfrac{G_{u-}}{(G^*_{u})^2}dG^*_u\;,
\end{array}$$
one can see that
$$\begin{array}{ll}Z^*_t
&=1-\dint_0^{t\wedge\ttau}\int_\bb R\dfrac{1}{G^*_s}p_s(s,l)\eta(ds,dl)\\
&\qquad-\dint_0^{t\wedge\ttau}\dfrac{1}{G^*_u}\Big\{\dint_0^u\dint_\bb
Rp_s(s,l)Z^{\theta_1,\theta_2}_{s,l}(u-)\theta_1(u;s,l)\eta(ds,dl)\Big\}'d
W(u)\\
&\qquad-\dint_0^{t\wedge\ttau}\dint_E\dfrac{1}{G^*_u}\Big\{\dint_0^u\dint_\bb
Rp_s(s,l)Z^{\theta_1,\theta_2}_{s,l}(u-) \theta_2(u,y;s,l))
\eta(ds,dl)\Big\}\{\mu(du,dy)-\nu(du,dy)\}\\
&\qquad-\dint_0^{t\wedge\ttau}\dfrac{G_{u-}}{(G^*_{u})^2}dG^*_u-\dfrac{G_{\ttau}}{G^*_\ttau}\bb
I_{t\ge\ttau}\\\\
&\qquad+
p_\ttau(\ttau,L)\bb I_{t\ge \ttau}+\dint_0^tp_{u-}(\ttau,L)\theta_1(u;\ttau,L)'\bb I_{u>\ttau}dW(u)\\
&\qquad +\dint_0^t\int_Ep_{u-}(\ttau,L)\theta_2(u,y;\ttau,L)\bb
I_{u>\ttau}\{\mu(du,dy)-\nu(du,dy)\}
\end{array}$$

$$\begin{array}{ll}
&=1-\dint_0^tZ^*_{u-}\bigg\{\dfrac{1}{G_{u-}}\dint_0^u\dint_\bb
Rp_s(s,l)Z^{\theta_1,\theta_2}_{s,l}(u-)\theta_1(u;s,l)\eta(ds,dl)\bb
I_{u\le\ttau}\\
&\qquad\qquad\qquad-\theta_1(u;\ttau,L)\bb I_{u>\ttau}\bigg\}'dW(u)\\
&\qquad-\dint_0^{t}\dint_EZ^*_{u-}\bigg\{\dfrac{1}{G_{u-}}\dint_0^u\dint_\bb
Rp_s(s,l)Z^{\theta_1,\theta_2}_{s,l}(u-) \theta_2(u,y;s,l))
\eta(ds,dl)\bb I_{u\le\ttau}\\
&\qquad\qquad\qquad -\theta_2(u,y;\ttau,L)\bb
I_{u>\ttau}\bigg\}\{\mu(du,dy)-\nu(du,dy)\}\\\\

&\qquad-\dint_0^{t\wedge\ttau}\dfrac{G_{u-}}{(G^*_{u})^2}dG^*_u-\dfrac{G_{\ttau}}{G^*_\ttau}\bb
I_{t\ge\ttau}\\
&\qquad+ p_\ttau(\ttau,L)\bb I_{t\ge
\ttau}-\dint_0^{t\wedge\ttau}\int_\bb
R\dfrac{1}{G^*_s}p_s(s,l)\eta(ds,dl),
\end{array}$$
which completes the proof.
\end{proof}

By Theorem \ref{prop-presentation of Z*}, we give the
$\bb G$-decomposition of a $(P,\bb F)$ martingale as follows:
\begin{theorem}\label{theorem-decomposition of F-martingale}
Let $p_s(s,l)$, $\theta_1(u;s,l)\bb I_{u>s}$ and
$\theta_2(u,y;s,l)\bb I_{u>s}$ be given as in Theorem
\ref{theorem-p(s;s,l)}. If $m$ is a c\`{a}dl\`{a}g $(P, \bb F)$-local
martingale of the following form
$$m_t=m_0+\int_0^t\xi_1(u)'dW(u)
+\int_0^t\int_E\xi_2(u,y)\{\mu(du,dy)-\nu(du,dy)\},$$
then
$$\begin{array}{ll}
X_t:=m_t+\dint_0^{t\wedge\ttau}\dfrac{1}{G_{u-}}\bigg\{\dint_0^u\dint_\bb
Rp_s(s,l)Z^{\theta_1,\theta_2}_{s,l}(u-)\Big\{\xi_1(u)'\theta_1(u;s,l)\\
\qquad\qquad\qquad\qquad\qquad\qquad\qquad\qquad+\dint_E\theta_2(u,y;s,l))\xi_2(u,y)F_u(dy)\Big\} \eta(ds,dl) \bigg\}du\\
\qquad\qquad-\dint_\ttau^t\bigg\{\xi_1(u)'\theta_1(u;\ttau,L)
+\dint_E\theta_2(u,y;\ttau,L)\xi_2(u,y)F_u(dy)\bigg\}du,
\end{array}$$
is a $(P,\bb G)$-local martingale.
\end{theorem}
\begin{remark}
\begin{enumerate}
  \item[1)]Theorem \ref{theorem-decomposition of F-martingale} may be viewed as
a corollary of
Callegaro, Jeanblanc and Zargari(2010)$^{\text{\cite{Callegaro-Jeanblanc-Zargari(2010)}}}$
and El Karoui, Jeanblanc and Jiao(2010)$^{\text{\cite{El
Karoui-Jeanblanc-Jiao2010}}}$, the main difference is that the
decomposition of $\bb F$-local martingale in $\bb G$ we give here
only depends on $p_s(s,l)$, $\theta_1(u;s,l)\bb I_{u>s}$ and
$\theta_2(u,y;s,l)\bb I_{u>s}$, since in the integral
$$\dint_0^u\dint_\bb
Rp_s(s,l)Z^{\theta_1,\theta_2}_{s,l}(u-)\Big\{\xi_1(u)'\theta_1(u;s,l)+\dint_E\theta_2(u,y;s,l))\xi_2(u,y)F_u(dy)\Big\}
\eta(ds,dl),$$ $\theta_1(u;s,l)=\theta_1(u;s,l)\bb I_{u>s}$ and
$\theta_1(u,y;s,l)=\theta_2(u,y;s,l)\bb I_{u>s}$, which is quite
interesting.
  \item[2)]Furthermore, since a $\bb G^{\ttau,L}$-stopping time might not be a
$\bb G$-stopping time and the optional projection of a $(P,\bb
G^{\ttau,L})$-local martingale on $\bb G$ might not be a $(P,\bb
G)$-local martingale(cf.Stricker(1977)).
Hence the proof of Proposition 3.3 in
Callegaro, Jeanblanc and Zargari(2010) is not strict. We can also use
the proof of Proposition 5.9 in  El Karoui, Jeanblanc and Jiao(2010) to
write $X_t$ into the following form $$X_t=X_1(t)\bb
I_{t<\ttau}+X_2(t;\ttau,L)\bb I_{t\ge\ttau}$$ where
$$\begin{array}{ll}
X_1(t)=m_t+\dint_0^{t}\dfrac{1}{G_{u-}}\bigg\{\dint_0^u\dint_\bb
Rp_s(s,l)Z^{\theta_1,\theta_2}_{s,l}(u-)\Big\{\xi_1(u)'\theta_1(u;s,l)\\
\qquad\qquad\qquad\qquad\qquad\qquad\qquad+\dint_E\theta_2(u,y;s,l))\xi_2(u,y)F_u(dy)\Big\} \eta(ds,dl)
\bigg\}du\bigg\}\\
X_2(t;s,l)=\{m_t-X_1(s)\}\bb
I_{t>s}-\dint_s^t\bigg\{\xi_1(u)'\theta_1(u;s,l)
+\dint_E\theta_2(u,y;s,l)\xi_2(u,y)F_u(dy)\bigg\}du,
\end{array}$$
and show that $\big\{X_2(t;s,l)p_t(s,l);t\ge s\big\}$ is a $(P,\bb
F)$-\textbf{local} martingale and $\big\{X_1(t)G_t +
\int_0^t\int_\bb R X_2(u;u, l)p_u(u, l)\eta(du,dl); t \ge
  0\big\}$
is a $(P,\bb F)$-\textbf{local} martingale. However, we can not use
Theorem \ref{theorem-G-martingale character}, since \textbf{the
sequence of the localization stopping times of
$\{X_2(t;s,l)p_t(s,l);t\ge s\}$ depends on $(s,l)$ which is
uncountable infinite}, one can see that the proof of  Proposition
5.9 in El Karoui-Jeanblanc-Jiao(2010) is not strict. Here we would
like to provide a strict proof based on Proposition
\ref{prop-presentation of Z*}.
\end{enumerate}
\end{remark}

\begin{proof}[Proof of Theorem \ref{theorem-decomposition of F-martingale}]
Let $m$ be a $(P,\bb F)$-local martingale, then $m$ is a $(P^*,\bb
F)$-local martingale which is also a $(P^*,\bb G^{\ttau,L})$-local
martingale.
Since
$$\begin{array}{ll}
\dfrac{1}{L^*_t}=\dfrac{G_t}{G^*_t}\bb I_{t<\ttau}+p_t(\ttau,L)\bb
I_{t\ge\ttau}=Z^*_t.
\end{array}$$
To prove $X$ is a $(P,\bb G)$-local martingale, we need only to
prove that $X_tZ^*_t$ is a $(P^*,\bb G)$-local martingale. As a
matter of fact, one can see from It\^{o}'s formula that
$$\begin{array}{ll}
X_tZ^*_t=m_0+\dint_0^tZ^*_{u-}dm_u+\dint_0^tZ^*_{u-}\bigg\{\dfrac{1}{G_{u-}}\dint_0^u\dint_\bb
Rp_s(s,l)Z^{\theta_1,\theta_2}_{s,l}(u-)\Big\{\xi_1(u)'\theta_1(u;s,l)\\
\qquad\qquad\qquad\qquad\qquad\qquad\qquad\qquad+\dint_E\theta_2(u,y;s,l))\xi_2(u,y)F_u(dy)\Big\} \eta(ds,dl)\bb
I_{u\le\ttau} \\
\qquad\qquad-\Big\{\xi_1(u)'\theta_1(u;\ttau,L)
+\dint_E\theta_2(u,y;\ttau,L)\xi_2(u,y)F_u(dy)\Big\}\bb
I_{u>\ttau}\bigg\}du\\
\qquad\qquad+\dint_0^tX_{u-}dZ^*_u+[X,Z^*]_t\\

=m_0+\dint_0^tZ^*_{u-}dm_u+\dint_0^tX_{u-}dZ^*_u\\
\qquad\qquad+\dint_0^tZ^*_{u-}\bigg\{\dfrac{1}{G_{u-}}\dint_0^u\dint_\bb
Rp_s(s,l)Z^{\theta_1,\theta_2}_{s,l}(u-)\Big\{\xi_1(u)'\theta_1(u;s,l)\\
\qquad\qquad\qquad\qquad\qquad\qquad\qquad\qquad+\dint_E\theta_2(u,y;s,l))\xi_2(u,y)F_u(dy)\Big\} \eta(ds,dl)\bb
I_{u\le\ttau} \\
\qquad\qquad\qquad\qquad-\Big\{\xi_1(u)'\theta_1(u;\ttau,L)
+\dint_E\theta_2(u,y;\ttau,L)\xi_2(u,y)F_u(dy)\Big\}\bb
I_{u>\ttau}\bigg\}du\\
\qquad\qquad-\dint_0^tZ^*_{u-}\bigg\{\dfrac{1}{G_{u-}}\dint_0^u\dint_\bb
Rp_s(s,l)Z^{\theta_1,\theta_2}_{s,l}(u-)\theta_1(u;s,l)\eta(ds,dl)\bb
I_{u\le\ttau}\\
\qquad\qquad\qquad-\theta_1(u;\ttau,L)\bb I_{u>\ttau}\bigg\}'\xi_1(u)du\\
\qquad\qquad-\dint_0^{t}\dint_EZ^*_{u-}\bigg\{\dfrac{1}{G_{u-}}\dint_0^u\dint_\bb
Rp_s(s,l)Z^{\theta_1,\theta_2}_{s,l}(u-) \theta_2(u,y;s,l))
\eta(ds,dl)\bb I_{u\le\ttau}\\
\qquad\qquad\qquad -\theta_2(u,y;\ttau,L)\bb
I_{u>\ttau}\bigg\}\xi_2(u,y)\mu(du,dy)
\end{array}$$

$$
\begin{array}{ll}
=m_0+\dint_0^tZ^*_{u-}dm_u+\dint_0^tX_{u-}dZ^*_u\\
\qquad\qquad-\dint_0^{t}\dint_EZ^*_{u-}\bigg\{\dfrac{1}{G_{u-}}\dint_0^u\dint_\bb
Rp_s(s,l)Z^{\theta_1,\theta_2}_{s,l}(u-) \theta_2(u,y;s,l))
\eta(ds,dl)\bb I_{u\le\ttau}\\
\qquad\qquad\qquad -\theta_2(u,y;\ttau,L)\bb
I_{u>\ttau}\bigg\}\xi_2(u,y)\{\mu(du,dy)-\nu(du,dy)\}.
\end{array}$$
Since both $m$ and $Z^*$ are $(P^*,\bb G)$-local martingale, one can
see that $XZ^*$ is a $(P^*,\bb G)$-local martingale, which completes
the proof.
\end{proof}

\begin{corollary}
Assume conditions of Theorem \ref{theorem-decomposition of F-martingale} hold, and let
$$\begin{array}{lr}
W^\bb
G(t)&:=W(t)+\dint_0^t\bigg\{\dfrac{1}{G_{u-}}\dint_0^u\dint_\bb
Rp_s(s,l)Z^{\theta_1,\theta_2}_{s,l}(u-)\theta_1(u;s,l)\eta(ds,dl)\bb
I_{u\le\ttau}\\
&-\theta_1(u;\ttau,L)\bb I_{u>\ttau} \bigg\}du,
\end{array}$$ then $W^\bb G$ is a $(P,\bb G)$-Brownian motion.
\end{corollary}

\begin{corollary}
Assume conditions of Theorem \ref{theorem-decomposition of F-martingale} hold, and let
$$\begin{array}{lr}
\nu^\bb G(du,dy)&:=\bigg\{1-\dfrac{1}{G_{u-}}\dint_0^u\dint_\bb
Rp_s(s,l)Z^{\theta_1,\theta_2}_{s,l}(u-)\theta_2(u,y;s,l)\eta(ds,dl)\bb
I_{u\le\ttau}\\
&+\theta_2(u,y;\ttau,L)\bb I_{u>\ttau} \bigg\}F_u(dy)du,
\end{array}$$ then $\nu^\bb G(du,dy)$ is the compensator of $\mu(du,dy)$ with respect to $(P,\bb G)$.
\end{corollary}
\begin{remark}
From Theorem \ref{theorem-decomposition of F-martingale},  we get the $\bb G$-decomposition of a $\bb F$ martingale, the
$(P,\bb G)$-Brownian motion and
 the compensator of $\mu(du,dy)$ with respect to $(P,\bb G)$  explicitly.
\end{remark}

\section{The weak representation of a $(P,\bb G)$-martingale}

Now, we have the following representation for a $(P^*,\bb G)$-martingale:

\begin{theorem} \label{theorem-present-(P^*,G) martingale}
Let $M_t:=M_1(t)\bb I_{t<\ttau}+M_2(t;\ttau,L)\bb I_{t\ge\ttau}$ be a
$(P^*,\bb G)$  martingale, then there exists a $\bb
G$-predictable process $\xi(u)$ with $\xi\in L^2_{loc}(W)$ and a $\wt{\scr P}(\bb G)$-measurable
function $\zeta(u,y)$ with $\zeta\in G_{loc}(\mu)$  such that
\begin{equation}\label{eqn-(P*,G) martingale-present}
\begin{array}{ll}M_t=M_0+\dint_0^t\xi(u)'dW(u)+\dint_0^t\int_E\zeta(u,y)\{\mu(du,dy)-\nu(du,dy)\}\\
\qquad+(M_2(\ttau;\ttau,L)-M_1(\ttau-))\bb
I_{t\ge\ttau}-\dint_0^{t\wedge\ttau}\dint_\bb R\dfrac{(M_2(u;u,l)-M_1(u-))}{G^*_u}\eta(du,dl),
\end{array}
\end{equation}
\end{theorem}
\begin{proof}
From Corollary \ref{theorem-(P*,G)-martingale character}, one can
see that  $\{M_1(t)G^*_t + \int_0^t\int_\bb R M_2(u;u,
l)\eta(du,dl); t \ge
  0\}$ and  $\{M_2(t;s, l); t\ge s\}$ are
$(P^*, \bb F)$-martingales for $\eta$-almost every $u\ge 0$ and $l\in
\bb R$, thus there exist  $\bb F$-predictable processes $\xi_1(u)$ and $\xi_2(u,s,l)$ with
 $\xi_1,\xi_2(;s,l)\in L^2_{loc}(W)$
 and $\wt{\scr P}(\bb F)$-measurable functions
$\zeta_1(u,y)$ and $\zeta_2(u,y;s,l)$ with $\zeta_1, \zeta_2(;s,l)\in G_{loc}(\mu)$ such that
$$\begin{array}{ll}
M_1(t)G^*_t + \dint_0^t\int_\bb R M_2(u;u,
l)\eta(du,dl)\\
\qquad\qquad=M_1(0)+\dint_0^t\xi_1(u)'dW(u)+\dint_0^t\dint_E\zeta_1(u,y)\{\mu(du,dy)-\nu(du,dy)\}\text{ and}\\
M_2(t;s,l)=M_2(s;s,l)+\dint_s^t\xi_2(u;s,l)'dW(u)+\dint_s^t\dint_E\zeta_1(u,y;s,l)\{\mu(du,dy)-\nu(du,dy)\}.
\end{array}$$
Thus
$$\begin{array}{ll}
M_1(t)G^*_t
=M_1(0)+\dint_0^t\xi_1(u)'dW(u)+\dint_0^t\dint_E\zeta_1(u,y)\{\mu(du,dy)-\nu(du,dy)\}\\
\qquad\qquad\qquad- \dint_0^t\int_\bb R M_2(u;u, l)\eta(du,dl),
\end{array}$$
one can see from It\^{o}'s formula that
$$\begin{array}{ll}
M_1(t)=M_1(t)G^*_t\dfrac{1}{G^*_t}\\
=M_1(0)+\dint_0^t\dfrac{\xi_1(u)'}{G^*_{u-}}dW(u)+\dint_0^t\dint_E\dfrac{\zeta_1(u,y)}{G^*_{u-}}\{\mu(du,dy)-\nu(du,dy)\}\\
\qquad\qquad- \dint_0^t\int_\bb R \dfrac{M_2(u;u,
l)}{G^*_{u-}}\eta(du,dl)+\dint_0^tM_1(u-)G^*_{u-}d\Big(\dfrac{1}{G^*_u}\Big)\\
=M_1(0)+\dint_0^t\dfrac{\xi_1(u)'}{G^*_{u-}}dW(u)+\dint_0^t\dint_E\dfrac{\zeta_1(u,y)}{G^*_{u-}}\{\mu(du,dy)-\nu(du,dy)\}\\
\qquad\qquad- \dint_0^t\int_\bb R \dfrac{M_2(u;u,
l)}{G^*_{u-}}\eta(du,dl)-\dint_0^t\dfrac{M_1(u-)}{G^*_u}dG^*_u\;.
\end{array}$$
Thus
$$\begin{array}{ll}
M_t=M_1(t)\bb I_{t<\ttau}+M_2(t;\ttau,L)\bb I_{t\ge\ttau}\\
=M_1(t\wedge\ttau)-M_1(\ttau)\bb I_{t\ge\ttau}+M_2(t;\ttau,L)\bb
I_{t\ge\ttau}\\
=M_1(0)+\dint_0^{t\wedge\ttau}\dfrac{\xi_1(u)'}{G^*_{u-}}dW(u)+\dint_0^{t\wedge\ttau}\dint_E\dfrac{\zeta_1(u,y)}{G^*_{u-}}\{\mu(du,dy)-\nu(du,dy)\}\\
\qquad\qquad- \dint_0^{t\wedge\ttau}\int_\bb R \dfrac{M_2(u;u,
l)}{G^*_{u-}}\eta(du,dl)-\dint_0^{t\wedge\ttau}\dfrac{M_1(u-)}{G^*_u}dG^*_u\\
\qquad\qquad-M_1(\ttau)\bb I_{t\ge\ttau}+M_2(t;\ttau,L)\bb
I_{t\ge\ttau}\\
=M_0+\dint_0^t\xi(u)'dW(u)+\dint_0^t\int_E\zeta(u,y)\{\mu(du,dy)-\nu(du,dy)\}\\
\quad\qquad+M_2(\ttau;\ttau,L)\bb
I_{t\ge\ttau}-\dint_0^{t\wedge\ttau}\dint_\bb R\dfrac{M_2(u;u,l)}{G^*_u}\eta(du,dl)\\
\quad\qquad-M_1(\ttau)\bb
I_{t\ge\ttau}-\dint_0^{t\wedge\ttau}\dfrac{M_1(u-)}{G^*_u}dG^*_u\\
=M_0+\dint_0^t\xi(u)'dW(u)+\dint_0^t\int_E\zeta(u,y)\{\mu(du,dy)-\nu(du,dy)\}\\
\quad\qquad+(M_2(\ttau;\ttau,L)-M_1(\ttau-))\bb
I_{t\ge\ttau}-\dint_0^{t\wedge\ttau}\dint_\bb R\dfrac{(M_2(u;u,l)-M_1(u-))}{G^*_u}\eta(du,dl)\;,
\end{array}$$
where
$$\begin{array}{ll}
\xi(u)=\dfrac{\xi_1(u)}{G^*_{u-}}\bb
I_{u\le\ttau}+\xi_2(u;\ttau,L)\bb I_{u>\ttau}\;,\\
\zeta(u,y)=\dfrac{\zeta_1(u,y)}{G^*_{u-}}\bb
I_{u\le\ttau}+\zeta_2(u,y;\ttau,L)\bb I_{u>\ttau}\;.
\end{array}$$
Since the $(P^*,\bb F)$-martingale $\{M_1(t)G^*_t + \int_0^t\int_\bb
R M_2(u;u, l)\eta(du,dl); t \ge
  0\}$ has no jumps at $\ttau$ as a $(P^*,\bb G)$-martingale and $G^*$
is continuous, one can see that $M_1(\ttau)=M_1(\ttau-)$, a.s. and
(\ref{eqn-(P*,G) martingale-present}) follows.
Note that since $(G_t)_t\geq 0$ is continuous
and never $0$, $\frac{1}{G}$ is locally bounded, then $\xi\in L^2_{loc}(W)$ and $\zeta\in G_{loc}(\mu)$, which complete the proof.
\end{proof}

Now we turn to  prove the predictable representation theorem  for a $(P,\bb G)$-martingale.
Similar to Lemma \ref{lemma-Af}, we have the following lemma
\begin{lemma} \label{lemma-Bf}
For any positive $\scr O(\bb F)\times\scr B$-measurable function
$f(s,l)$ such that $E[\dint_{0}^{\infty}\int_\bb Rf(s,l)$
$p_s(s, l)\eta(ds,dl)]<\infty$, let
$$A^f_t:=\int_0^t\int_\bb R \dfrac{f(s,l)}{G_{s-}}p_s(s,l)\eta(ds,dl),$$
then $A^f$ is a continuous increasing process and
$$M^f_t=f(\ttau,L)\bb I_{t\ge\ttau}-A^f_{t\wedge\ttau}$$
is a $(P,\bb G)$-martingale, i.e., $A^f_{\0\wedge\ttau}$ is the
$(P,\bb G)$-compensator of $f(\ttau,L)\bb I_{t\ge\ttau}$.
\end{lemma}
\begin{proof} Similar to the proof of Lemma \ref{lemma-Af}, one can
see that for any $t_1<t_2$,
$$\begin{array}{ll}
E\big[f(\ttau,L)\bb I_{t_1<\ttau\le t_2}\big|\scr F_{t_1}\big]
&=\dint_{t_1}^{t_2}\int_\bb Rf(s,l)p_{t_1}(s,l)\eta(ds,dl)\\
&=E\Big[\dint_{t_1}^{t_2}\int_\bb
Rf(s,l)p_{s}(s,l)\eta(ds,dl)\Big|\scr F_{t_1}\Big]\end{array}$$ and
that
$$\begin{array}{ll}
E\Big[\dint_{t_1}^{t_2}\int_\bb R\dfrac{\bb
I_{s\le\ttau}}{G_{s-}}p_s(s,l)f(s,l)
\eta(ds,dl)\Big|\scr F_{t_1}\Big]\\
=E\Big[\dint_{t_1}^{t_2}\int_\bb R\dfrac{E[\bb I_{s\le\ttau}|\scr
F_s]}{G_{s-}}p_s(s,l)f(s,l)
\eta(ds,dl)\Big|\scr F_{t_1}\Big]\\
=\dint_{t_1}^{t_2}\int_\bb Rf(s,l)p_s(s,l)\eta(ds,dl),
\end{array}$$
and the rest is completely the same as the proof of Lemma
\ref{lemma-Af}.
\end{proof}

Similarly to the process of the prove of Theorem \ref{theorem-present-(P^*,G) martingale} , we can prove  the following theorem.
\begin{theorem}\label{theorem-present-(P,G) martingale}
Let $M_t:=M_1(t)\bb I_{t< \ttau}+M_2(t;\ttau,L)\bb I_{t\geq\ttau}$
be a  $(P,\bb G)$-martingale, then there exists a $\bb G$-predictable
process $\xi(u)$ with $\xi\in L^2_{loc}(W)$ and a $\wt{\scr P}(\bb G)$-measurable function
$\zeta(u, y)$ with $\zeta\in G_{loc}(\mu)$ such that
\begin{equation}\label{M-01}
\begin{array}{ll}
M_t=M_0+\dint_0^t\xi(u)'dW^\bb G(u)+\dint_0^t\int_E\zeta(u,y)\{\mu(du,dy)-\nu^\bb G(du,dy)\}\\
\qquad+(M_2(\ttau;\ttau,L)-M_1(\ttau-))\bb
I_{t\ge\ttau}-\dint_0^{t\wedge\ttau}\dint_\bb R\dfrac{(M_2(u;u,l)-M_1(u-))}{G_{u-}}p_u(u,l)\eta(du,dl).
\end{array}
\end{equation}
\end{theorem}

\section{Conclusion}
The achievement of such a study will provide
practitioners with specific formulas to deal with the progressively enlarged filtration issued
from a Brownian motion and an integer-valued random measure. In this paper, the formulas
are presented as functionals of three parameters
$(p_s(s,l), \theta_1, \theta_2)$.
In this paper,
we deeply characterize the conditional density process
and give the Doob-Meyer's decomposition of the survival process.
We also discuss the  necessary and sufficient conditions
for a  $\bb G$-martingale.
By  the  Lemma \ref{lemma-Af}, we explicitly have described
the $\bb G$-decomposition of a $(P, \bb
F)$-martingale and  have proved the martingale representation theorems.   Formula parametrization in the enlarged filtration is a useful quality in financial modeling and
we will consider more applications in the next work.

\vskip15pt

\end{document}